\input amstex
\input psfig.sty
\documentstyle{amsppt}
\topmatter
\title
Injectivity Radius and Fundamental Groups of Hyperbolic 3-Manifolds
\endtitle
\rightheadtext{Injectivity Radius and Fundamental Groups}
\author
Matthew E. White
\endauthor
\abstract
It is shown that for each integer $n>1$ there exists a constant
$R_n>0$ such that if $M$ is a closed hyperbolic 3-manifold with 
Rank $\pi_1(M) = n$, then the injectivity radius of $M$ is bounded above
by $R_n$.
\endabstract
\date
26 October 1998
\enddate
\endtopmatter
\document
\head 1. Introduction
\endhead

Thurston's Geometrization conjecture implies that every closed irreducible
3-manifold with infinite fundamental group containing no ${\Bbb Z}\oplus
 {\Bbb Z}$
subgroup is hyperbolic.  Thus, it is of interest to obtain information about the
hyperbolic metric from purely topological data.  For example, the Gromov norm,
which is defined homologically, turns out to be a certain scalar multiple of 
volume if the manifold is hyperbolic [Th].  Similarly, there is also an upper
bound on the volume in terms of length of any presentation of its fundamental
group [C]. The {\bf injectivity radius} of a closed hyperbolic 3-manifold $M$,
denoted $inj(M)$, is one half the minimum length of all essential loops in $M$.
This is equivalent to the usual definition from differential geometry (see
[DC]). Recall that the {\bf rank} of a finitely generated group is the minimum
number of elements required to generate the group.  In this paper we give an
upper bound on the injectivity radius in terms of fundamental group rank.  It
is clear that no such lower bound exists since Dehn-filling of a hyperbolic
manifold provides examples of bounded rank with injectivity radius approaching
zero.  In a subsequent paper, we hope to provide a lower bound on injectivity
radius in terms of other group theoretic data.  The basic outline of the proof
in our present setting is the following:  given a closed hyperbolic 3-manifold
$M$, we construct a graph of minimal length which carries $\pi_1(M)$.  By
considering the graph's preimage in hyperbolic space, we show that if the
injectivity radius is sufficiently large, any relation in the fundamental group
forces some of the graph's edges into close proximity.  This in turn allows
modification of the graph to reduce its length while still carrying fundamental
group.  We therefore show that, in fact, $\pi_1(M)$ is a free group of finite
rank.  A closed hyperbolic 3-manifold is irreducible. Therefore, its
fundamental group cannot be a free product.  This contradiction proves:
\proclaim{Theorem 4.4} For each integer $n>1$ there exists a universal constant
$R_n>0$ such that if $M$ is a closed hyperbolic 3-manifold and Rank $\pi_1(M) =
n$, then $inj(M) < R_n$. \endproclaim \noindent Residual finiteness of closed
hyperbolic 3-manifolds then provides us with: \proclaim{Corollary 4.5} Given a
closed hyperbolic 3-manifold $M$ and an integer $n>0$, there exists a finite
sheeted cover $$p:\tilde M\longrightarrow\ M$$ with $Rank\ \pi_1 (\tilde M) \ge
n$. \endproclaim

We shall use $M$ to denote a closed, connected hyperbolic 3-manifold.
In particular, $M$ is isometric to a quotient of ${\Bbb H}^3$ by a group
of loxodromic isometries.  We say that a connected finite graph $\Gamma$
is an {\bf $n$-graph} if all vertices of $\Gamma$ are trivalent and
$Rank\ \pi_1(\Gamma) = n$.  A {\bf carrier $n$-graph} for a manifold $M$
is an $n$-graph $\Gamma$ together with a map $f:\Gamma\longrightarrow M$
such that $f_*:\pi_1(\Gamma)\longrightarrow \pi_1(M )$ is an
epimorphism.  In this case, the {\bf length} of $\Gamma$ in $M $ is
defined by
$$\ell(f(\Gamma)) = \sum_{e\ \text{an edge of}\ \Gamma} |e| $$
where $|e|$ is the length of $e$ measured by pulling back the path
metric of $M $ to $\Gamma$. When we refer to a
carrier $n$-graph for a given manifold $M$, we
implicitly assume that $n = Rank\ \pi_1(M)$.

This paper is organized as follows: In section 2, we prove some
necessary technical results using trivalent graphs.  Section 3 contains
the $n=2$ case of the main theorem; this was previously known to C.
Adams.  The main theorem and  its corollary are
proved in section 4.  We remark here that our
results extend naturally to the case of bounded
negative curvature.

\head 2. Trivalent Graphs of Minimal Length
\endhead
\medskip
\proclaim{Proposition 2.0}
For each integer $n>1$, the set of $n$-graphs is finite.  Each $n$-graph
has $3(n-1)$ edges.
\endproclaim
\medskip
\noindent {\bf Proof} If $\Gamma$ is an $n$-graph, then $$\chi(\Gamma) = 1-n =
(\text{number of vertices}) - (\text{number of edges})$$.  Note that twice the
number of edges equals three times the number of vertices since all vertices of
$\Gamma$ are trivalent.  Therefore, there are $3(n-1)$ edges. $\square$
\medskip
\noindent Let $e_1\dots e_{3n-1}$ be the oriented edges of $\Gamma$.
Given a loop in $\Gamma$, its homotopy class may be represented by a {\bf
directed edge path}. This is a loop in $\Gamma$ which is a product of
homeomorphisms $\gamma_i :[0,1]\longrightarrow e_i$ with the property that for
each $i\ge 2$, we have $\gamma_{i-1}(1)=\gamma_i(0)$.  In our case, we shall
also assume that the maps $\gamma_i$ are parameterized proportionally to arc
length.  Although a slight abuse of notation, the standard convention is to
denote a directed edge path as $$(e_{i_1}^{\pm 1},e_{i_2}^{\pm
1},\dots,e_{i_k}^{\pm 1})$$ with $e_i^{-1}$ denoting $e_i$ with the opposite
orientation.  A directed edge path that is a closed loop is called {\bf
reduced} if $e_{i_m} \ne e_{i_{m-1}}^{-1}$ for any index $i$ and $e_{i_1}\ne
e_{n_k}^{-1}$.  For further details, see [M].

\medskip
\proclaim{Proposition 2.1} Let $f:\Gamma\longrightarrow M $ be a carrier $n$-graph.  If $\alpha$ is a simple closed curve in
$\Gamma$, then $f_*($[$\alpha$]$) \ne 0$.  In particular,
$\ell(f(\alpha)) > inj(M)$.
\endproclaim
\medskip\noindent {\bf Proof}  Suppose $\alpha$ is a simple closed curve
in $\Gamma$.  Then $\alpha$ represents a generator of $\pi_1(\Gamma)$.
Therefore, $f_*($[$\alpha$]$) \ne 0$ as otherwise we could generate
$\pi_1(M)$ with fewer than $n$ elements.  Moreover, the length of any
essential loop in $M $ is at least $2 inj(M)$. $\square$

\proclaim{Lemma 2.2}
A closed hyperbolic 3-manifold $M$ has a carrier $n$-graph $\Gamma$ of
minimal length.  The edges of $\Gamma$ map to geodesic segments in $M$.
\endproclaim
\medskip\noindent {\bf Proof}
We shall omit much detail since the proof closely follows the standard proof
of the existence of closed geodesics in free homotopy classes (see [DC]).  Fix
$M$ and let $\Gamma$ be an $n$-graph.  We may regard $\Gamma$ as a compact
subset of $\Bbb R^k$ with each edge smoothly embedded.  Put
$$d = \text{inf }\lbrace\ell (f(\Gamma))\ \vert \
f:\Gamma\longrightarrow M \ ,\ f_*\text{ is an epimorphism}\rbrace.$$
\noindent 
Now suppose that $f:\Gamma\longrightarrow M$ is a carrier $n$-graph.  Suppose
also that there is a loop in $\Gamma$ containing exactly one edge.  By
Proposition 2.1, the restriction of $f$ to this edge must lift to a path in
${\Bbb H}^3$ with distinct endpoints.  Therefore, if the restriction of $f$ to
an arbitrary edge of $\Gamma$ is not constant, it must lift to a path in ${\Bbb
H}^3$ with distinct endpoints.  This means we may construct a family of maps
$F=\lbrace f_i:\Gamma\longrightarrow M \rbrace_{i=1}^\infty$ such that for each
$i\ge 1$, $f_{i*}$ is an epimorphism, $f_i$ is geodesic on each edge of
$\Gamma$, and $\ell(f_i(\Gamma)) \longrightarrow d$. Since $M$ is compact, a
straightforward Arzela-Ascoli argument then shows that the closure of $F$ is
compact in ${\Cal C}[\Gamma,M]$.  Thus, there exists a subsequence ${f_j}$ of
$F$ which converges uniformly to a continuous map $f:\Gamma\longrightarrow M$.
Uniform convergence implies that $f$ maps every edge to a geodesic segment and
that $\ell(f(\Gamma)) = d$.  We show that $f_*$ is an epimorphism. Choose a
vertex $x_0\in\Gamma$.  Given $\epsilon$ small enough, cover each edge of
$f(\Gamma)$ with finitely many $\epsilon$ balls such that each is isometric to
a ball in ${\Bbb H}^3$.  Since the convergence of the $f_j$ is uniform,  by
picking $j$ large enough, we can ensure that $f_j$ maps into this covering.
Now $f_{j*}$ is an epimorphism.  Thus, given $[\beta]\in\pi_1(M,f_j(x_0))$, let
$[\alpha]\in\pi_1(\Gamma,x_0)$ be such that $f_{j*}([\alpha])=[\beta]$.  Note
that every homotopy class in $\pi_1(\Gamma,x_0)$ is represented by a directed
edge path in $\Gamma$.  Hence, we assume that $\alpha$ is a directed edge path.
In each ball, we may homotop $f_j$ to $f$.  This means $f(\alpha)$ is freely
homotopic to $f_j(\alpha)$.  In particular, we have that $f_*$ is an
epimorphism.  By proposition 2.0, there are finitely many $n$-graphs.  Hence we
many repeat the above proof for each distinct $n$-graph and then select the
graph of minimal length.$\square$ \medskip Minimal length carrier graphs have a
very nice symmetry property which we shall later put to good use. \medskip
\proclaim{Definition 2.3} A subset of a hyperbolic 3-manifold is a {\bf
Y-subset} if it is isometric to a subset ${\Cal S}$ of ${\Bbb H}^3$ with the
following properties.\newline \noindent (i) ${\Cal S}$ is composed of three
geodesic arcs which meet only at a common endpoint.\newline \noindent (ii) Each
pair of geodesic arcs in ${\Cal S}$ meets the common endpoint with incidence
angle $120^\circ$. \endproclaim \noindent It follows at once that a Y-subset is
planar. \medskip \proclaim{Lemma 2.4} Given a manifold $M $ with Rank
$\pi_1(M)\ge 2$, there exists a minimal length carrier $n$-graph
$f:\Gamma\longrightarrow M $ with the following property: each vertex $v$ of
$\Gamma$ has a neighborhood N($v$) such that $f$ maps $N(v)$ isometrically to a
Y-subset of $M$. In particular, $f(v)$ is the common endpoint of the geodesic
edges in $f(N(v))$. Furthermore, every edge has non-zero length.\endproclaim
\medskip \noindent {\bf Proof} We assume that $f:\Gamma\longrightarrow M $ is a
fixed minimal length carrier $n$-graph with the property that $f$ is an
arc-length parameterization on each edge of $\Gamma$.  Let $v$ be a vertex of
$\Gamma$.  Suppose first that no edge of $\Gamma$ has length 0 under $f$.
Choose a small neighborhood $W$ of $f(v)$ such that $W$ is isometric to an open
ball in ${\Bbb H}^3$.  The edges of $\Gamma$ incident at $v$ map to geodesic
edges in $M$.  If any pair of these edges intersect in $W$ at a point other
than $f(v)$, then the image of one edge is contained in that of the other edge.
It is easy to see that in this case we can modify the graph to reduce its
length. Therefore, the images of the three edges intersect only at $f(v)$ in
$W$.  It also follows that these three edges are planar inside $W$, since
otherwise orthogonal projection to a hyperbolic plane reduces the length of
$\Gamma$.  We now suppose there exists a pair of edges with incidence angle
less than $120^\circ$.  We modify the map $f$ near $v$ as pictured in Figure
2.0.  A short hyperbolic trigonometry calculation then shows that this reduces
total length.

\medskip
\psfig{file=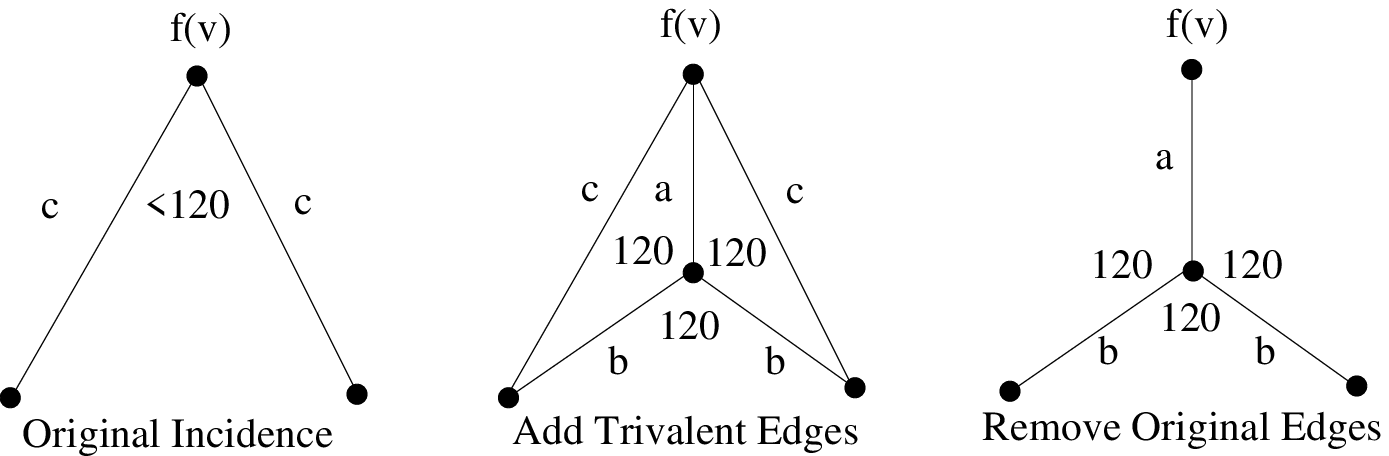,height=1.5in,width=5.0in}
Figure 2.0. Adjustment near a vertex: $2c > a + 2b$.
\medskip

\noindent Thus, every incidence angle is at least $120^\circ$.  Of course, the
trivalent angle sum is exactly $360^\circ$.   To complete the proof in the
general case, we show that if there exists an edge $e$ of $\Gamma$ with length
0, then we may modify $f$ to reduce total length. Suppose such an edge $e$
exists. Form a maximal connected subgraph $T$ such that $T$ contains $e$ and
$\ell(f(T)) = 0$.  We note that Proposition 2.1 implies that $T$ must be a
tree. Hence there exist two vertices $v_1$ and $v_2$ in $T$ such that
$f(v_1)=f(v_2)$ and each vertex has exactly one incident edge of length $0$.
\medskip
\psfig{file=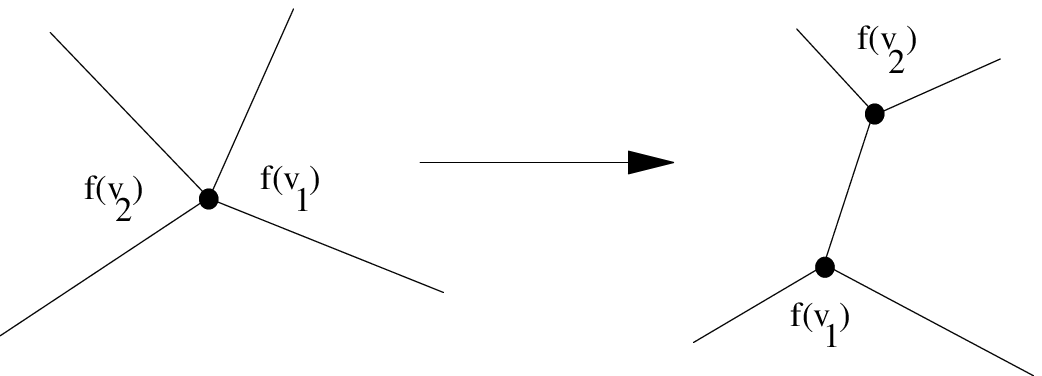,height=1.5in,width=5.0in}
Figure 2.1.  Adjustment at identified vertex.
\medskip

\noindent There are two possibilities.  In case one, there are four edges
incident at $f(v_1)=f(v_2)$.  Then we can modify the image of $f$ using the
trivalent structure exactly as shown in Figure 2.1. This move separates the
vertices by introducing a new edge.  The idea is to stretch the edge
of length $0$ to have small positive length. Then we modify the map using the
move in Figure 2.0. This produces a new map of a (possibly different) carrier
$n$-graph of strictly lower length.  In case two, one of the edges incident at
$f(v_2)$ is a subset of an edge incident at $f(v_1)$ and all other incidence
angles are $120^\circ$.  In this case, a short calculation shows that we may
reduce length by sliding $f(v_1)$ along the identified edges while
straightening the remaining two. $\square$ \medskip

To simplify the exposition, we note the following fact.  A minimal length carrier
$n$-graph need not be embedded; however, intersections are transverse.  Thus,
given any $\epsilon > 0$ we can adjust the graph by pushing apart intersecting
edges inside an $\epsilon$ ball.  This allows us to produce an embedded graph
with arbitrarily small total increase in length.  Of course, some edges will
now fail to be geodesic within finitely many $\epsilon$ neighborhoods. But the
distortion can be made as small as desired. It is our aim to show that given an
integer $n>0$, there is a constant $R_n$ such that if a manifold $M $ has
$inj(M ) > R_n$, then the minimal length carrier $n$-graph
$f:\Gamma\longrightarrow M$ provided in Lemma 2.2 has $f_*$ an isomorphism.
This will imply that $\pi_1(M)$ is a free group of rank $n$ which contradicts
that $M$ is a closed hyperbolic 3-manifold.  To do this, we must prove that
$f_*$ is injective.  Hence, we try to understand the generic properties of a
non-trivial class $[\alpha]\in \pi_1(\Gamma)$ for which $f_*([\alpha]) = 0$.
\proclaim{Definition 2.5} Fix a universal covering $p:{\Bbb H}^3\longrightarrow
M$ and let $f:\Gamma\longrightarrow M$ be a minimal length carrier $n$-graph.
\newline (a) We say that $\gamma :S^1\longrightarrow \Gamma$ is a {\bf
compressing loop} if $\gamma$ is essential in $\Gamma$ and $f\gamma$ is null
homotopic in $M $.\newline (b) A {\bf shortest compressing loop} is a
compressing loop of shortest length in $\Gamma.$\newline (c) A {\bf standard
lift} for $kernel(f_*)$ is a lift to $\Bbb H^3$ of $f\gamma :S^1\longrightarrow
M $ where $\gamma :S^1\longrightarrow \Gamma$ is a shortest compressing loop.
\endproclaim \medskip

\noindent {\bf Remarks} It follows easily that we may represent a
shortest compressing loop by a reduced directed edge path.  Thus, Lemma
2.4 and the remarks above show we may assume that a standard lift for
$kernel(f_*)$ is a piecewise geodesic embedding of $S^1$ in $\Bbb H^3$
in which the geodesic pieces meet at $120^\circ$ angles.  We shall call
the geodesic pieces in the image of a standard lift {\bf edges}.
Notice that these are lifted edges of the minimal length carrier graph.
We now abstract the important properties of standard lifts.

\proclaim {Definition 2.6}
A path $h:[0,1]\longrightarrow {\Bbb H}^3$ is a {\bf Geodesic-120 path}
if the following conditions hold: \newline
(i)  $h$ is a (possibly closed) piecewise geodesic path. \newline
(ii) edges in the image of $h$ meet at $120^\circ$ angles. \newline
\endproclaim
\medskip
The key point is that any standard lift we construct for a given closed
hyperbolic 3-manifold will be a closed geodesic-120 path in $\Bbb H^3$.
For our purposes, it will be helpful to think of the image of a standard
lift as a tractable geometric object in ${\Bbb H}^3$.  Frequently, we
shall also need to consider an important type of subpath in a
geodesic-120 path.
\medskip
\proclaim {Definition 2.7}
Let $h:[0,1]\longrightarrow {\Bbb H}^3$ be a geodesic-120 path.  Let
$[a,b]\subset [0,1]$ where $h(a)$ and $h(b)$ are endpoints of edges
in $H = h([0,1])$.  A {\bf segment} of $h$ joining $h(a)$ to
$h(b)$ is a subpath $s=h|_{[a,b]}:[a,b]\longrightarrow {\Bbb H}^3$.
\endproclaim
\medskip
Naturally, if a geodesic-120 path is an embedding, given two endpoints
there is exactly one segment joining them.
\proclaim{Lemma 2.8}
Let $f:\Gamma\longrightarrow M $ be a minimal length carrier $n$-graph.
If $h:S^1\longrightarrow\Bbb H^3$ is a standard lift, then in any
segment of $h(S^1)$ containing at least $3(n-1)$ edges,
there is a edge with hyperbolic length $\ge {inj(M)\over {3(n-1)}}$.
\endproclaim
\medskip
\noindent {\bf Proof} Given a manifold $M $, Proposition 2.0 shows that
an $n$-graph $\Gamma$ in $M $ has $3(n-1)$ edges.
Thus, a simple closed curve in $\Gamma$ has less than $3(n-1)$ edges.
By Proposition 2.1, the length of such a curve is at least $inj(M)$. 
Now since $h$ is a standard lift, $h=f\gamma$ where $\gamma$ is a shortest
compressing loop.  We consider how $\gamma$ behaves in $\Gamma$. We may view
$\gamma$ as a product of directed closed edge paths
$\alpha_1\alpha_2\ldots\alpha_k$ with the following properties: (a) each path
$\alpha_i$ has at most $3(n-1)$ edges (b) each path $\alpha_i$ contains an edge
$e$ such that the path defined by removing the edge $e$ from $\alpha_i$ is not
a loop (c) each path $\alpha_i$ contains a simple closed curve
with at most $3(n-1)$ edges. Thus, each  $\alpha_i$ contains an edge of
length at least $inj(M )\over {3(n-1)}$.  The proof is complete since this
argument is independent of the starting vertex for the path $\gamma$. $\square$

\bigskip
\head 3. The Case $Rank(\pi_1(M))=2$.
\endhead
In the remainder of this paper, we use the Poincare Disc model for $\Bbb
H^3$. We shall also let  ${\Cal W}$  denote the horoball with diameter
given by the geodesic joining 0 to (-1,0,0). In this section, we shall
prove the main theorem in the special case of Rank $\pi_1(M)$ = 2. Note
that there are precisely two 2-graphs (see figure 3.0).
\bigskip

\psfig{file=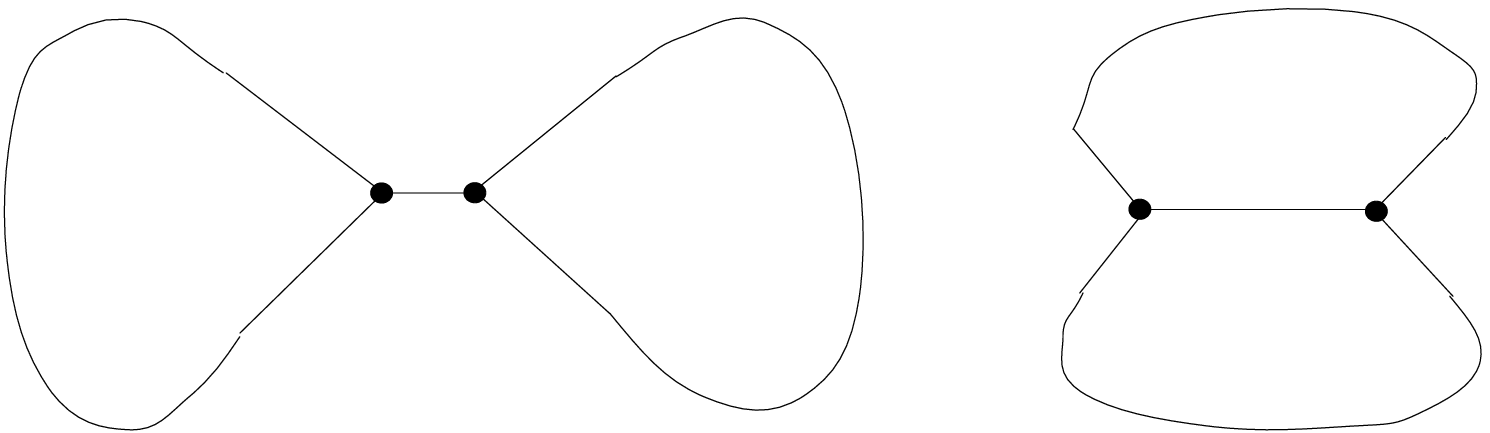,height=1.5in,width=5.0in}
Figure 3.0. The two possible 2-graphs.
\medskip

\proclaim {Lemma 3.0} There exists $R_2>0$ such that if $M$ is a closed,
connected, hyperbolic 3-manifold with Rank $\pi_1(M)$ = 2 then $inj(M)
\le R_2$.
\endproclaim
\medskip
\noindent {\bf Proof}  The proof is established by considering
Geodesic-120 paths in $\Bbb H^3$.  Given a closed, connected hyperbolic
3-manifold $M$ and a minimal length carrier $2$-graph
$f:\Gamma_2\longrightarrow M$ we may choose a standard lift
$h:S^1\longrightarrow {\Bbb H}^3$ for $kernel (f_*)$.  Let  $H$ denote
$h(S^1)$.  Since $H$ is compact, there exist points $v$,$w \in H$ of
maximal distance apart.  By an isometry of $\Bbb H^3$, we may assume that
$w=0$ and that $v$ lies on the geodesic joining 0 to (-1,0,0).  Now
$H\subset \bar B(v,d(v,0))$ which implies that  $H$  is contained in the
horoball  ${\Cal W}$ .  Note also that since  $H$  is a Geodesic-120
path, 0 and $v$ must be endpoints of edges of  $H$ .  Moreover, since
$H$  is a loop, there must be two edges of  $H$ incident at 0.
Consider these two edges.  Geodesics through 0 in the disc model
correspond to Euclidean diameters of the unit ball in $\Bbb R^3$.  Since
the angle between the edges is $120^\circ$,  there exists $L_0>0$ so that if
one of the edges of $H$  incident at 0 has length at least $L_0$, then  $H$  is
not contained in  ${\Cal W}$.  There are precisely two 2-graphs (see figure
3.0) and each of these graphs has three edges. By using Proposition 2.1, we see
that a 2-graph may have at most one ``short''  edge.  More precisely, if
$inj(M) > 2L_0$, then the image $H$ of any standard lift has the following
property: given any pair of edges of $H$ with a common endpoint, at least one
edge has length greater than $L_0$.  In other words, put $R_2=2L_0$.   If
$inj(M) > R_2$, (after isometry of ${\Bbb H}^3$) any standard lift for
$kernel(f_*)$ must have a edge which is not contained in  ${\Cal W}$.  This
contradiction shows that $kernel(f_*)$ is trivial, so that $\pi_1(M)$ is
actually a free group of rank 2.  Since this is impossible, we conclude that
$inj(M) \le R_2$. $\square$

\medskip
\psfig{file=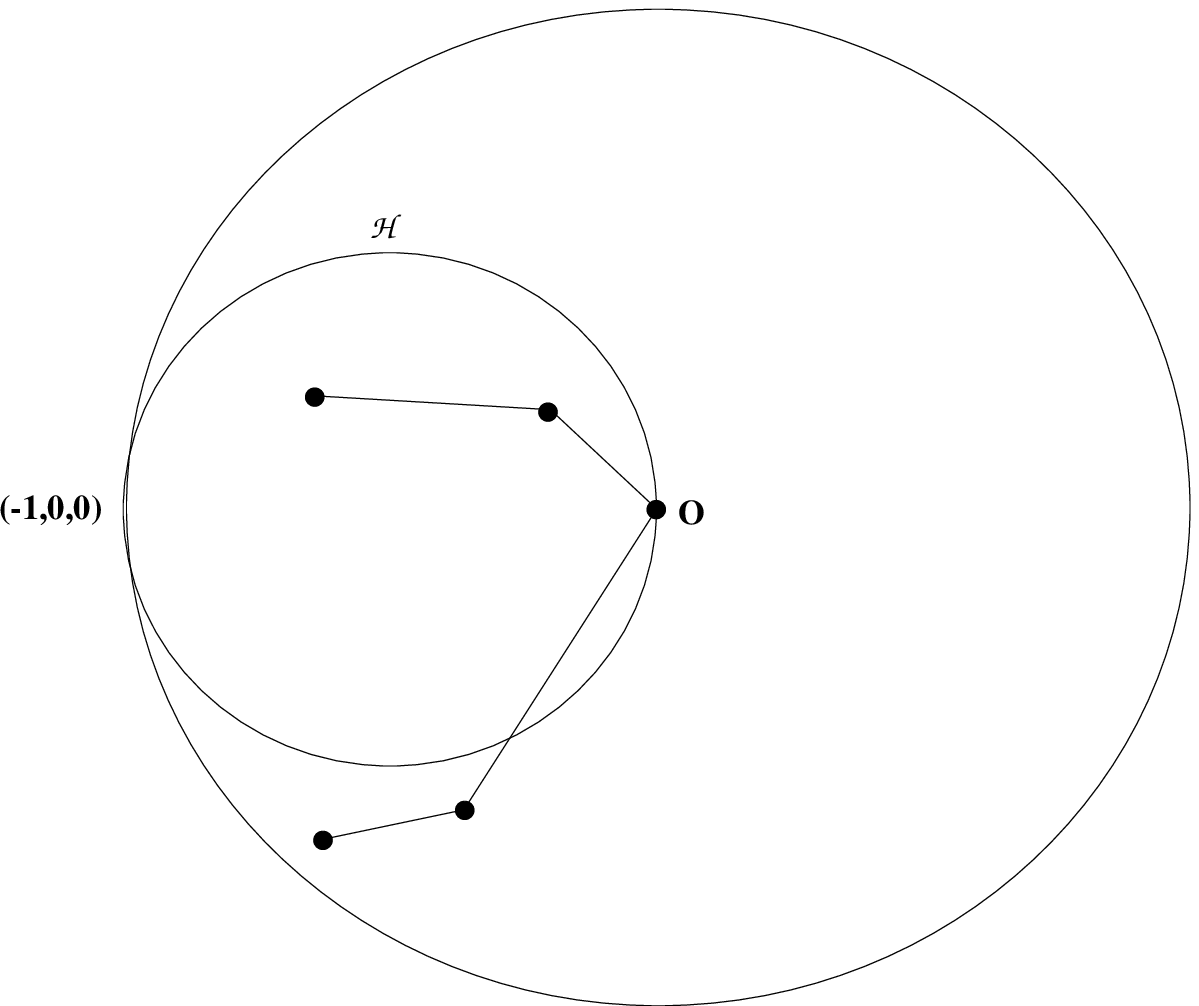,height=4.0in,width=5.0in}
Figure 3.1 Escape of Standard Lift in Rank 2 Case.
\medskip

\medskip\noindent {\bf Remark}  It is possible to explicitly compute an upper
bound for injectivity radius in the $n=2$ case.  One finds using a very
simple geometry calculation that $L = 2.6338\ldots$ so that $inj(M) \le
5.2676\ldots$ \bigskip

\head 4. Proof of Main Theorem
\endhead
The general case is more subtle because $n$-graphs with $n>2$ have more
edges.  Although large injectivity radius will ensure that minimal
length carrier graphs have some long edges, if $n>2$ these graphs may
have subpaths of very short edges. This prevents us from using the
simple geometric approach in Lemma 3.0 since the image of a standard
lift can have very large geodesic edges while remaining inside the
horoball ${\Cal W}$ .  Instead, we will argue that, for sufficiently
large injectivity radius, if the standard lift is contained in  ${\Cal
H}$, then it is possible to reduce the length of the minimal length
carrier $n$-graph.  We will produce a {\bf short cut arc} in ${\Bbb
H}^3$ which allows us to make a shorter carrier graph using
``cut-and-paste'' in the manifold below.  To do this, we need two
technical lemmas which we now describe.  Recall that ${\Bbb H}^3$ enjoys
``thin triangles.''  This means there is a universal constant $\Delta$, the
{\bf thin triangles constant}, such that given a hyperbolic triangle with
geodesic sides $X$, $Y$ and $Z$ and a point $x\in X$, there exists a point
$y\in Y\cup Z$ with $d(x,y)<\Delta.$  The idea is to use this fact together
with the geometry of the horoball to show that sufficiently large injectivity
radius forces a standard lift to be ``thin.''   Our first lemma is motivated by
this fact:

\medskip
\proclaim {Proposition 4.0} Given $L>0$, if $M$ is a closed hyperbolic
3-manifold with $n =$ Rank $(\pi_1(M))$, and $inj(M) > [3(n-1)]^2 L$, then
every standard lift for $\pi_1(M)$ has two edges of length at least $L$.
\endproclaim \medskip\noindent {\bf Proof} Let $M$ satisfy the hypothesis, and
let $h:S^1\longrightarrow\Bbb H^3$ be a standard lift with image $H$.  Lemma
2.8 shows that any standard lift has one edge of length at least $3(n-1)L$. But
if $H$ has exactly one edge of length at least $3(n-1) L$, then it has at most
$3(n-1)$ edges.  Therefore, $H$ cannot be a closed loop in ${\Bbb H}^3$ unless
there is another edge of length at least $L$. $\square$ \medskip \noindent This
proposition suggests we should try to produce our short cut arc between two
``long'' edges.  The lemma below says that, if $A$ and $B$ are two such edges
in ${\Cal W}$ that start at points ``close'' to the origin, then there is a
short cut arc between them.

\medskip
\proclaim {Lemma 4.1 (Short Cut Lemma)}
For any $\delta>0$ there exists $\bar L(\delta)>0$ such that if $A=[a$,$a^\prime]$
and $B=[b$,$b^\prime]$ are geodesic edges contained in  ${\Cal W}$
with:\newline
(i)  $a$,$b\in \bar B(0,\delta),$\newline
(ii) $|A|$,$|B| \ge \bar L$,\newline
\noindent then there exist points $e\in A$ and $f\in B$ such
that:\newline
(i)  $d(e$,$a) < {1 \over 3}\bar L\ \text{and}\ d(f$,$b) < {1 \over
3}\bar L$\newline
(ii) $d(e$,$f) < d(e$,$a)-\Delta$ \text{and}\ $d(e$,$f) < d(f$,$b)-\Delta$.
\endproclaim
\medskip\noindent {\bf Proof}  See Figure 4.0. Given $A$ and $B$ as in the
hypothesis, let $l_A=[0$,$a^\prime]$, $l_B=[0$,$b^\prime]$, and $\phi$ be the
angle between $l_A$ and $l_B$.  Let $\delta >0$.  Choose $L_1 >
2(\delta+\Delta)$.  This choice assures that given $A,B$ as in hypothesis
with $|A|$,$|B| > L_1$, we have $a^\prime$,$b^\prime \not \in \bar
B(0$,$\delta)$.  Now notice that as $L_1\rightarrow\infty$, we have 
for all $A,B$ as in the hypothesis with $|A|$,$|B| > L_1$, that 
$a^\prime$,$b^\prime\rightarrow$ (-1,0,0) in the Euclidean metric.
Hence, we may choose $L_2>L_1$ so that given any $A$ and $B$ as in the
hypothesis with $|A|$,$|B| > L_2$, we are guaranteed to have
$|l_A|$,$|l_B|$ sufficiently large and $\phi$ sufficiently small that
there exist points $x\in l_A$ and $y\in l_B$ with $d(x,\partial \bar
B(0$,$\delta))$, $d(y,\partial \bar B(0$,$\delta)) = 5\Delta$ and
$d(x$,$y)<\Delta$.  Using thin triangles, there exist points $e\in A$
and $f\in B$ such that $d(e$,$x)$,$d(f$,$y) < \Delta$.  This implies
$d(e$,$f) < 3\Delta$.  Now choose $\bar L(\delta) = 3max(L_2,6\Delta + 2\delta)$.
Then $d(e$,$a)$,$d(f$,$b) > 4\Delta$ which proves that
$d(e$,$f)<d(e$,$a)-\Delta$ and $d(e$,$f)< d(f$,$b)-\Delta$.  Also, $d(e$,$a)\le
6\Delta + 2\delta$ so that $d(e$,$a)< {1\over 3}\bar L(\delta)$. $\square$
\medskip
\noindent {\bf Remark} In subsequent work, we shall cut out one of the
geodesics $[e,a]$ or $[f,b]$ and paste in $[e,f]$.  It is worth mentioning at
this point that conclusion $(ii)$ tells us $\ell ([e,a])-\ell ([e,f])>\Delta$
and $\ell ([f,b])-\ell ([e,f])>\Delta$, so the short cut reduces length by at
least $\Delta$.  The small perturbation needed to make $\Gamma$ embedded can be
chosen to change lengths less than $\Delta$.

\medskip
\psfig{file=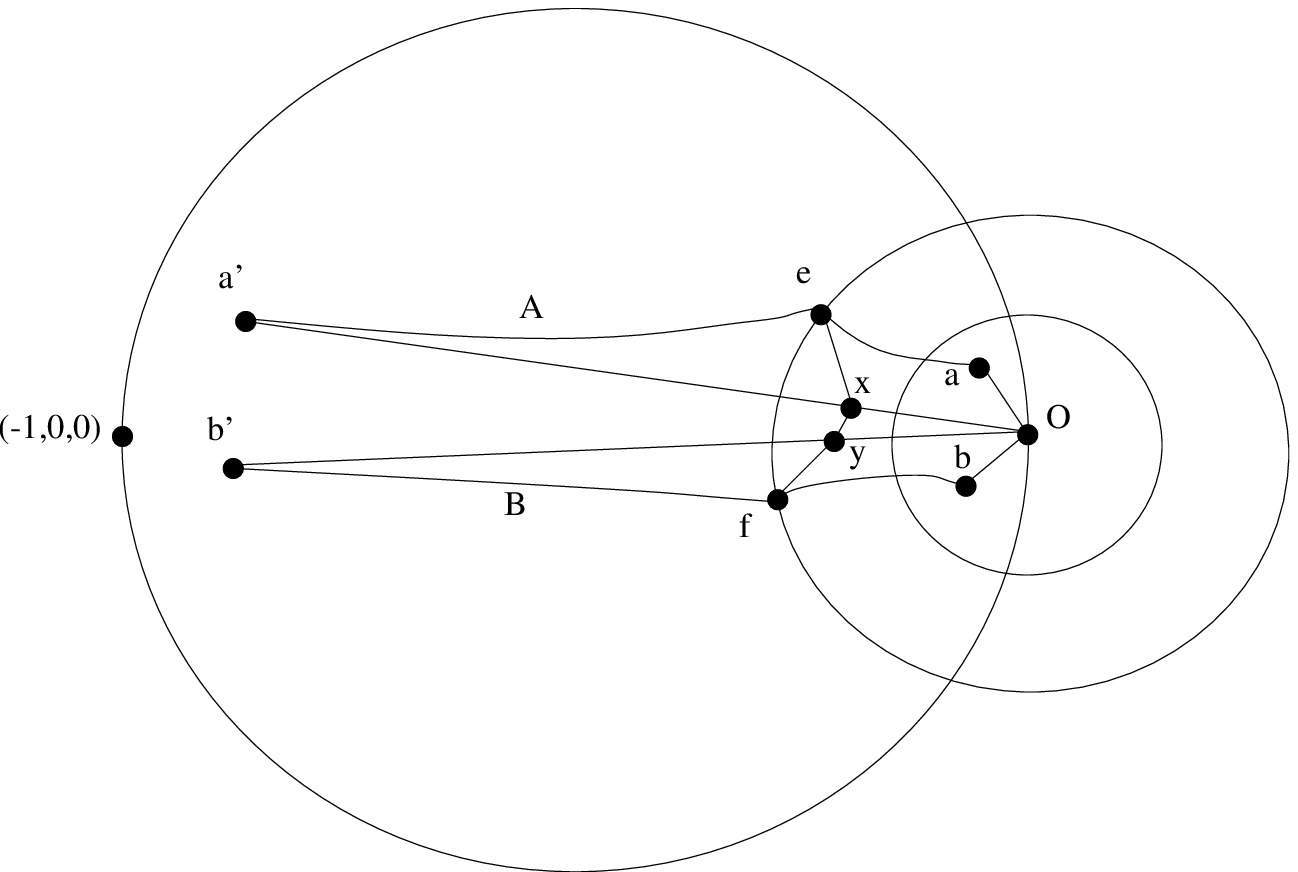,height=3.5in,width=5.0in}
Figure 4.0. Trapping of geodesic segments.
\medskip

\medskip The next two results show that we can actually achieve the conditions
in the hypothesis of Lemma 4.1  for embedded Geodesic-120 paths provided we
attach the right assumptions. Suppose that $n=Rank(\pi_1(M))$. If $inj(M)$ is
very large, Lemma 2.8 and Proposition 4.0 imply that there are two ``very
long'' edges in a standard lift for $\pi_1(M)$.  Denoting these two edges $A$
and $B$, we then notice that by Lemma 2.8 there is a Geodesic-120 path of at
most $6n-6$ edges joining $A$ to $B$.  We show that either the starting points
of $A$ and $B$ are sufficiently close or there is another (not quite as)
``long'' edge between $A$ and $B$.  Taking innermost ``long'' edges then
provides the requirement for the short cut lemma.  The basis for the procedure
is the following fact about finite sequences of numbers: \medskip \proclaim
{Proposition 4.2} Let $L(0)<L(1)<L(2)<\dots$ be a given sequence of positive
numbers.  Let $x_0,\dots x_n$ and $y_0,\dots y_m$ be sequences of positive
numbers with $k=n+m$ such that:\newline (i) $x_n$ and $y_m$ are both greater
than $max\lbrace x_{n-1},...,x_0,y{m-1},...,y_0\rbrace$. \newline (ii) $x_n$
and $y_m$ are both greater than $L(k)$\newline Then one of the following
conditions must hold:\newline (1) There exist $r\ge 0$ and $s\ge 0$ such that
for $q = r+s$, \newline \indent (i) $x_r$ and $y_s$ are both greater than
$max\lbrace x_{r-1},...,x_0,y{s-1},...,y_0\rbrace$. \indent (ii) $x_r$ and
$y_s$ are both greater than $L(q)$\newline \indent (iii) $\sum_{i=0}^{r-1} x_i
+ \sum_{i=0}^{s-1} y_i \le q L(q-1)$ \newline
(2) Either $x_0 > L(0)$ or $y_0 > L(0)$.
\endproclaim
\medskip
{\bf Proof} Let $x_0,\dots x_n$ and $y_0,\dots y_m$ be sequences which
satisfy the hypothesis.  We may assume that $n\ge 1$ and $m\ge 1$ as
otherwise $(2)$ holds.  We may also assume that $\sum_{i=1}^n x_i +
\sum_{i=1}^m y_i > k L(k-1)$.
Then without loss of generality there exists $x_r$ with $r< n$ such
that $x_r > L(k-1)$.  By choosing $s=m$ if necessary, there also exists
$y_s>L(k-1)$.  We choose $r$ and $s$ to be as small as possible.  If
$r=0$ or $s=0$, the proof is finished since $(2)$ holds.  Otherwise we
have that $x_r$ and $y_s$ are both greater than $max\lbrace
x_{r-1},...,x_0,y{s-1},...,y_0\rbrace$. Note also that $r+s\le n-1+m = k-1$.
Hence, for $q=r+s$ we have $x_r>L(k-1)\ge L(q)$ and $y_s>L(k-1)\ge
L(q)$.  Now if $1(iii)$ holds, the proof is done.  Suppose otherwise.
Repeat the above procedure on the sequences $x_r,\dots x_n$ and
$y_s,\dots y_m$; continuing in this manner, we either obtain sequences
$x_{r'},\dots x_0$ and $y_{s'},\dots y_0$ which satisfy $(1)$ or we run
out of sequence elements so that we obtain $(2)$. $\square$
\medskip
\proclaim {Lemma 4.3}
For each integer $k\ge 0$ there is an $L(k)>0$ such that if
$g:[0,1]\longrightarrow {\Bbb H}^3$ is an embedded Geodesic-120 path
with $G=g([0,1])$ which satisfies:\newline
(i) $0\in G$ \newline
(ii) There exist two geodesic edges $A$ and $B$ in $G$ with
$|A|$,$|B|>L(k)$\newline
(iii) Let $S$ be the segment of $g$ which joins $A$ to $B$.
Then $0\in S$ and $S$ has at most $k$ edges. \newline
(iv) Every edge of $S$ has length at most $L(k)$.
\smallskip\noindent Then one of the following holds:\newline
\noindent (1) $G$ is not contained in the horoball  ${\Cal W}$ .
\newline
\noindent (2) There exist edges $E_1$ and $E_2$ in $G$ such
that:\newline
(i)   $g^{-1}(E_1)$ and $g^{-1}(E_2)$ are in different components of
$[0,1]-{g^{-1}(0)}$\newline
(ii)  Every edge in the segment joining $E_1$ to $0$ has length
less than $|E_1|$ and $|E_2|$ \newline
(iii) Every edge in the segment joining $E_2$ to $0$ has length
less than $|E_1|$ and $|E_2|$ \newline
(iv)  Let $e_1\in E_1$ and $e_2\in E_2$ be the endpoints of these
edges contained in the segment joining $E_1$ to $E_2$.  Then
there is a short cut arc $F$ joining $e \in E_1$ to $f\in E_2$ such that
$d(e_1$,$e)>|F|+\Delta$, $d(e_1$,$e) < {1\over 3} L(k)$, and $d(e_2$,$f) <
{1\over 3} L(k)$.
\endproclaim

{\bf Proof} Define the following sequence:  set $L(0)$ equal to the $L_0$ in
the proof of Lemma 3.0.  Recall that this is the length required to force
one of the edges at $0$ to leave the horoball ${\Cal W}$. Put $L(1)=2L(0)$ and
for $k\ge 2$ define $L(k)=\bar L(kL(k-1))$ with $\bar L$ provided by Lemma 4.1.
This implies that for $k\ge 2$, $L(k)\ge 2(k L(k-1)+\Delta)$
so that evidently $L(0)<L(1)<L(2) ...$ is an increasing sequence of positive
numbers. Now fix $k\ge 0$ and let $G$ be the image of a geodesic-120 path which
satisfies the hypothesis.  Let $x_n = |A|$, $y_m = |B|$, and let $x_{n-1},\dots
x_0$ denote the lengths of the edges joining $A$ to $0$ labelled in decreasing
order as we proceed from $A$ to $0$. Likewise, let $y_{m-1},\dots y_0$ denote
the lengths of the edges joining $B$ to $0$.  By hypothesis $(iii)$, $n+m =
k_1\le k$.  By $(ii)$, $x_n>L(k)\ge L(k_1)$ and $y_m>L(k)\ge L(k_1)$.  Also, by
$(iv)$ we have that $x_n$ and $y_m$ are both greater than $max\lbrace
x_{n-1},...,x_0,y{m-1},...,y_0\rbrace$. This means these sequences satisfy the
hypothesis of Proposition 4.2.  Conclusion $(2)$ of Proposition 4.2 gives that
$x_0>L(0)$ or $y_0>L(0)$. This means that either $x_0$ or $y_0$ leaves the
horoball ${\Cal W}$; hence, $G$ is not contained in ${\Cal W}$. Otherwise,
conclusion $(1)$ of Proposition 4.2 gives us subsequences $x_r,\dots x_0$ and
$x_s,\dots x_0$ with these properties: $x_r>L(q)$, $y_s>L(q)$, $x_r$ and $y_s$
are both greater than $max\lbrace x_{r-1},...,x_0,y{s-1},...,y_0\rbrace$, and
$\sum_{i=0}^{r-1} x_i + \sum_{i=0}^{s-1} y_i \le q L(q-1)$.  Refer to Figure
4.1. Geometrically, this means the following.  Let $E_1$ and $E_2$ be the edges
of $G$ which correspond respectively to the lengths $x_r$ and $y_s$. It follows
at once that statements $2(i)$ through $2(iii)$ hold. To see that $2(iv)$ is
also true, note that the segment joining $E_1$ to $0$ has length at most
$qL(q-1)$ since $\sum_{i=0}^{r-1} x_i \le q L(q-1)$.  The analogous statement
holds for the segment joining $E_2$ to $0$.  Thus, $e_1,e_2\in \bar
B(0,qL(q-1))$.  Since $|E_1| = x_r > L(q)$ and $|E_2| = y_s > L(q)$, Lemma 4.1
applies to give us the required short cut arc $F$. $\square$

\bigskip
\psfig{file=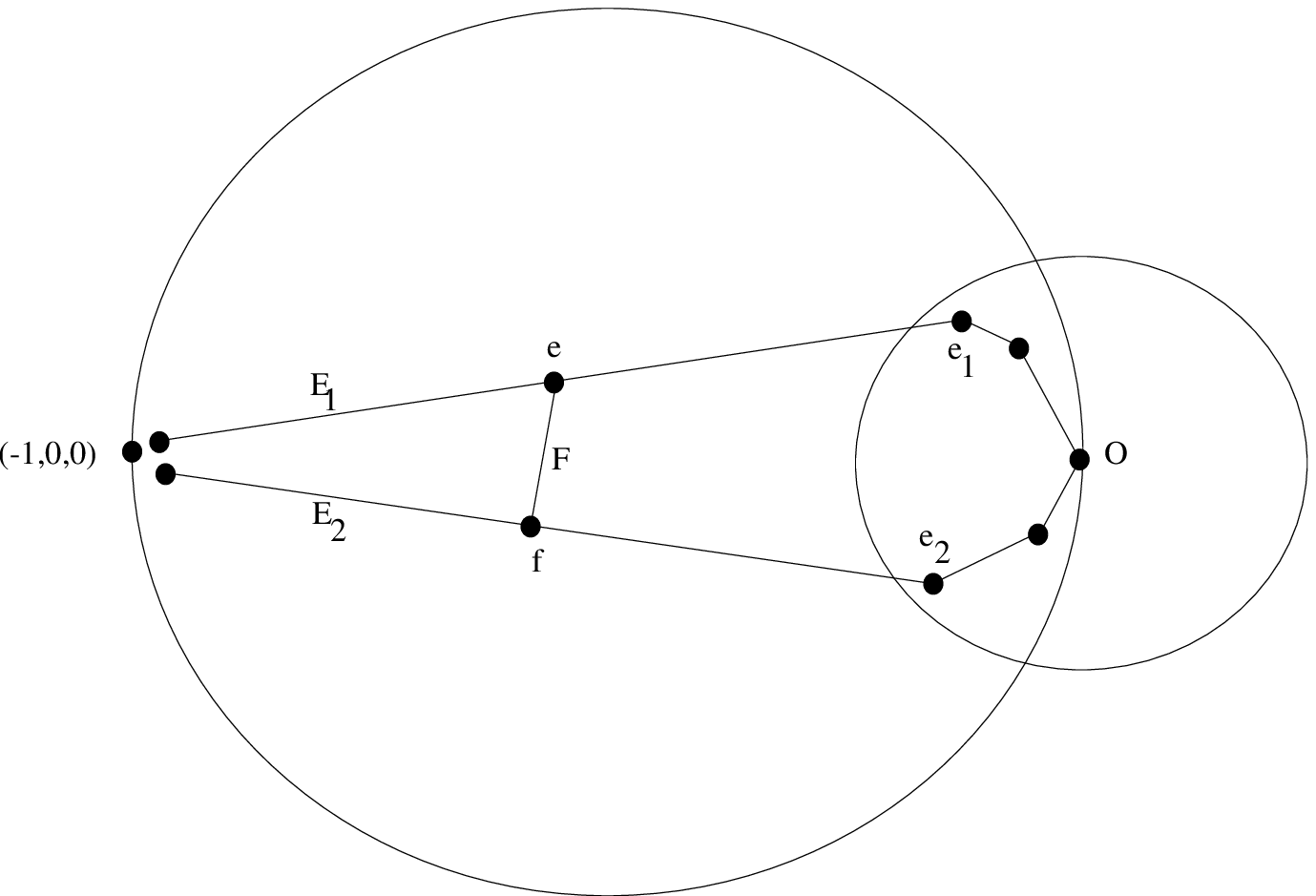,height=3.5in,width=5.0in}
Figure 4.1.
\medskip

We now show how to use Lemma 4.3 to establish the main result in this
paper. By the remark following the proof of Lemma 4.1, we can assume our
graph is embedded without affecting the validity of our ``cut-and
paste'' argument described below.  This means we can think of a standard
lift as a geodesic-120 embedding of $S^1$.
\medskip
\proclaim {Theorem 4.4}
For each integer $n>1$ there exists an $R_n > 0$ such that if M is a
closed, connected, hyperbolic 3-manifold and Rank $\pi_1$(M) = n, then
$inj(M) < R_n$.
\endproclaim
\medskip
\noindent {\bf Proof of Main Theorem} Given an integer $n>1$, we will
show that there exists $R_n>0$ such that if a closed hyperbolic
3-manifold $M $ has Rank $\pi_1(M )=n$ and $inj(M ) \ge R_n$, then
$\pi_1(M)$ is a free group.  This will establish the theorem since the
fundamental group of such a manifold cannot be a free group.   Fix
$n>0$.  We shall use the notation $L(k)$ as in the statement of Lemma
4.3. Notice that we may choose $R>0$ so that if $n=Rank(\pi_1(M))$ and 
$inj(M)>R$ then the following holds: if $f:\Gamma\longrightarrow M $ is 
a minimal length carrier graph, then every standard lift for $kernel(f_*)$
has two edges of length
at least $L=L(3(n-1))$.   Set $R_n$ = max(R, (3n-3)L(3n-4)).  
Thus,  suppose $M$ has $inj(M)>R_n$ and Rank $\pi_1(M )=n$.   Let $h$ be a
standard lift for $kernel(f_*)$ to $\Bbb H^3$.  As in the proof of lemma 3.0,
$h$ is a Geodesic-120 path through $0$ with $H = h(S^1)$ contained in the
horoball ${\Cal W}$.  Also, using the remark following Definition 2.5, we
assume $h$ is an embedding of $S^1$. Notice that there are two edges $A$ and
$B$ in $H$  with $|A|$,$|B|> L(3n-3)$ and at most $3n-3$ edges in the segment
joining them through $0$.  Consider the Geodesic-120 path $g$ defined by
$A$,$B$ and the segment of $h$ containing $0$ which joins them.  By the
proof of Lemma 3.0, since $|A|$ and $|B|$ are greater than $L(1)$, $g$
is not a closed loop.  Thus, $g$ satisfies the hypothesis of Lemma 4.3.
Therefore, the assumption that $H$ is contained in  ${\Cal W}$  means
there exist edges $E_1$,$E_2$ in $G = g([0,1])$ which satisfy $2.(i)
- 2.(iv)$ in Lemma 4.3.  Let $e_1$,$e_2$ denote the respective endpoints
of $E_1$,$E_2$ contained in the segment of $g$ which joins
$E_1$ to $E_2$.  This segment contains $0$.  Moreover, since
$H$  is contained in ${\Cal W}$, the proof of Lemma 4.3 shows that we may
assume that $E_1$ and $E_2$ are innermost with respect to $0$ in the
following sense: There is an integer $i>0$ with $|E_1|$,$|E_2| > L(i)$
and $e_1$,$e_2\in\bar B(0,iL(i-1))$ where $iL(i-1)<(3n-3)L(3n-4)$.  To
see this, note that if no such $i$ exists, then we may repeat the
inductive argument in the proof of Proposition 4.2 to obtain a edge
$E$ incident at $0$ with $|E|>L(0)$.  This contradicts the assumption
that  $H\subset  {\Cal W}$ .

\medskip
\psfig{file=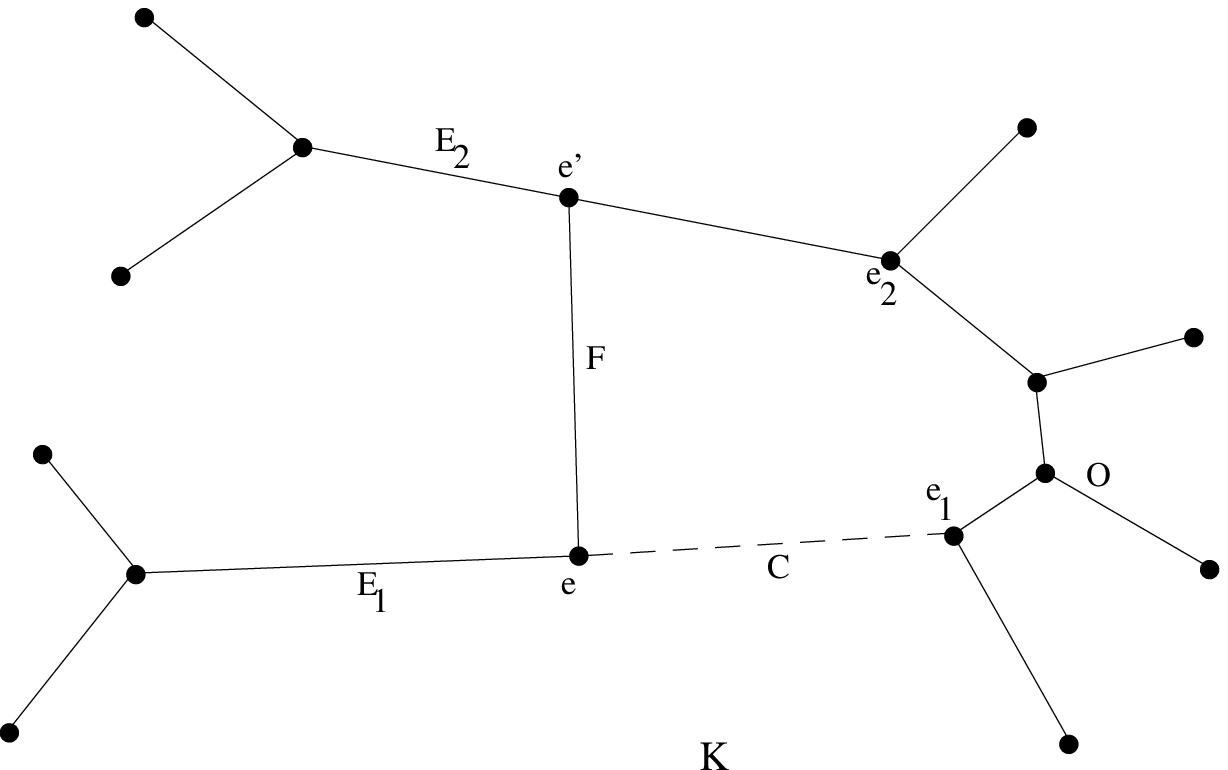,height=2.75in,width=5.0in}
Figure 4.2.
\medskip

\noindent Refer to figure 4.2.  We complete the proof using the arc $F$
joining $E_1$ to $E_2$ provided by Lemma 4.3.  We let $e\in E_1$ and
$e^\prime\in E_2$ denote the endpoints of $F$.  Let $C$ denote the
subarc of $E_1$ which joins $e_1$ to $e$. We wish to show that we can remove
$C$ and attach $F$ such that the projection to $M$ gives a shorter carrier
n-graph.  Let $\rho : {\Bbb H}^3\longrightarrow M$ denote the universal cover
of $M$.  We first show that the modified image $(f(\Gamma) - \rho (int(C)))
\cup \rho (F)$ defines a map of an $n$-graph into $M$.  We can work with the
image $f(\Gamma)$ since $\Gamma$ is embedded in $M$. Now, $\rho (e_1)$ is a
vertex of $f(\Gamma)$.  By removing $\rho (int(C))$, we convert $\rho(e_1)$
into a bivalent vertex. Hence, we can amalgamate the remaining two edges
incident at $\rho (e_1)$ into a single edge. This removes the vertex $\rho
(e_1)$.  We then view $\rho ((E_1 - int(C)) \cup F$ as a single edge so that
$\rho (e^\prime)$ becomes a trivalent vertex.  This new subset of $M$ is an
embedded, closed, trivalent graph $f^\prime : \Gamma^\prime\longrightarrow M$
with the same number of vertices as $\Gamma$; therefore, $f^\prime :
\Gamma^\prime\longrightarrow M$ is an $n$-graph. Note also that $\ell (f^\prime
(\Gamma^\prime)) < \ell (f(\Gamma))$.  Therefore, we need only show that
$\Gamma^\prime$ carries $\pi_1 (M)$. Consider again the carrier $n$-graph
$f:\Gamma\longrightarrow M $.  Let $P$ denote the path defined by the edges
joining $e_2$ to $e_1$ together with the subarc of $E_2$ joining $e'$ to $e_2$.
The idea is then to show that $P$ can be used to replace $C$ in any loop of
$\Gamma$. This follows if we show we have not cut edges in $P$ by
removing $C$.  Choose a basepoint $x_0\in f(\Gamma) - (\rho(E_1)\cup\rho(E_2))$
and let $y_0\in f^{-1}(x_0)$. Let $[\alpha]\in\pi_1(M,x_0)$.  Since $f_*$ is an
epimorphism, there exists $[\beta]\in\pi_1(\Gamma,y_0)$ with
$f_*([\beta])=[\alpha]$. We shall construct a corresponding loop $\beta^\prime$
in $\Gamma^\prime$ so that $f_*^\prime[\beta^\prime] = [\alpha]$. It follows
from $2.(i)$ through $2.(iv)$ of Lemma 4.3 that every edge in the segment of
$g$ which joins $E_1$ to $E_2$ has length less than both $E_1$ and $E_2$.
Therefore, none of these edges projects to $\rho(E_1)$ or to $\rho(E_2)$.  Now
suppose $\rho(E_1)=\rho(E_2)$. We must have $\rho(e_1) \ne \rho(e_2)$ since
otherwise $d(e_1,e_2) < (3n-3)L(3n-4)$ implies the geodesic edge with endpoints
$e_1$ and $e_2$ projects to a nontrivial loop in $M$ of length less than $R_n$.
This guarantees that the arc of $E_2$ joining $e'$ to $e_2$ does not project to
$\rho(C)$. Thus, $\rho(C)$ is not a subset of $\rho(P)$.  In particular,
$\rho(F\cup P)$ and $\rho(C)$ are paths in $M$ that are homotopic with
endpoints fixed.  These facts imply that if we imagine $\beta$ as a directed
edge path in $\Gamma$, we can use the word for $\beta$ to build $\beta^\prime$
in $\Gamma^\prime$.  That is, $\rho(E_1)$ corresponds to an edge $W_1$ in
$\Gamma$ and $\rho((E-1-int(C)) \cup F \cup P)$ corresponds to an edge path
$W_2$ in $\Gamma^\prime$.  In the word for $\beta$, replace $W_1$ (resp.
$W_1^{-1}$) by $W_2$ (resp. $W_2^{-1}$). This gives a loop $\beta^\prime$ in
$\Gamma^\prime$.  By construction, we may homotop $f\beta$ to
$f^\prime\beta^\prime$ fixing $x_0$, which shows that $f^\prime_*$ is onto.
This completes the proof of the main theorem. $\square$

\proclaim{Corollary 4.5} Given a closed hyperbolic 3-manifold $M$ and an
integer $n>0$, there exists a finite sheeted cover $$p:\tilde M\longrightarrow\
M$$ with $Rank\ \pi_1 (\tilde M) \ge n$. \endproclaim \medskip \noindent {\bf
Proof}.  Given $M$ and $n>0$, let $\tilde M_0$ be a covering space of $M$.
Suppose we have $Rank\ \pi_1 (\tilde M_0) < n$. Then Theorem 4.3 shows
$inj(\tilde M_0)< R_n$. $\tilde M_0$ has finitely many closed geodesics
$\gamma_1,\dots,\gamma_k$ of length less than $R_n$.  By residual finiteness,
there exists a a covering space  $\tilde M$ of $\tilde M_0$ such that the loops
$\gamma_1,\dots,\gamma_k$ do not lift.  This implies $inj(\tilde M) > R_n$ so
that $Rank\ \pi_1 (\tilde M_0) \ge n$. $\square$ \bigskip \head References
\endhead

\bigskip
\noindent {\bf Acknowledgements} The author thanks Professor Daryl Cooper for
his encouragement and for many helpful conversations.
\bigskip

\noindent [C] D. Cooper, ``The Volume of a Closed Hyperbolic
3-Manifold is Bounded by $\pi$ Times the Length of any Presentation of its
Fundamental Group'', To appear P.A.M.s. \smallskip

\noindent [DC] M. Do Carmo, ``Riemannian Geometry'', Birkhauser, Boston,
1993.
\smallskip
\noindent [E] P. Eberlein,  ``Structure of Manifolds of Nonpositive
Curvature'' in {\it Global Geometry and Global Analysis 1984, Proceedings
1984}, D. Ferus et.al. editors, Lecture Notes in Math. 1156, Springer, New
York, 1985, pp. 86-153. 1993. \smallskip \noindent [G] M. L. Gromov,
``Hyperbolic Groups'' in {\it Essays in Group Theory}, S. Gersten, MSRI
Publications 8, Springer, New York, 1987, pp. 75-263. 1993. \smallskip
\noindent [M] W. Massey, ``Algebraic Topology: An Introduction'',
Harbrace, New York, 1967.
\smallskip
\noindent [Pr] A. Preissman, {\it Quelques proprietes globales des espaces de
Riemann}, Comm. Math. Helv. 15 (1943), pp. 175-216. \smallskip
\noindent [Th] W.P. Thurston, ``The Geometry and Topology of
Three-Manifolds'', Princeton University, 1979.
\enddocument
\end